\documentclass[leqno,12pt]{article}

\usepackage{amsfonts,amssymb,amsmath,amsgen,theorem}

\global
\global
\theoremstyle{break} \theorembodyfont{\it}

\newtheorem{theoreme}{Theorem} \theorembodyfont{\sl}
 \theorembodyfont{\sl}

\newtheorem{proposition}{Proposition}

 \newtheorem{rem}{Remark}
\theoremheaderfont{\sc}

\newcommand{\R}{  \mathbb{R}   }

\begin{document}
 \bibliographystyle{plain}
\title{\vskip -1cm On uniqueness for the critical wave equation }

  \author{Nader Masmoudi\\
{\small Courant Institute of Mathematical Sciences,}\\
 {\small 251 Mercer Street, New York NY 10012}\\
{\small masmoudi@cims.nyu.edu}\\
and \\
Fabrice Planchon \\
{\small  Laboratoire Analyse, G\'eom\'etrie \& Applications}\\
{\small UMR 7539, Institut Galil\'ee}\\
{\small Universit\'e Paris 13, 99 avenue J.B. Cl\'ement}\\
{\small 93430 Villetaneuse FRANCE}\\
{\small fab@math.univ-paris13.fr}
}

 \date{}

 \maketitle
 \begin{abstract}
We prove the   uniqueness of 
  weak solutions  to the critical defocusing wave equation
  in 3D under a local energy inequality condition.
 More precisely, we prove the  uniqueness of    $  u  \in  L^\infty_t(\dot{H}^{1})\cap 
\dot{W}^{1,\infty}_t(L^2)$, under the condition  that $u$  verifies  some  local
  energy inequalities.
 \end{abstract}
 \par \noindent

\section{Introduction and statement of result}
\label{sec:intro}

We consider the defocusing quintic wave equation in 
3D, 
\begin{equation} \label{eq:5d}
\left\{  \begin{array}{l} 
  \Box u  + u^5=0, \\
   u(t=0) = u_0,   \quad u_t(t=0) = u_1.  
 \end{array}  \right. 
\end{equation}
Existence of  global weak solutions goes back to Segal (\cite{segal}, under milder
assumptions on the nonlinearity). Existence of global smooth solutions was proved by Grillakis
(\cite{Gri}), while global solutions in the energy space $C({\mathbb R};
H^1) \cap C^1({\mathbb R}; L^2) $ were constructed by Shatah and Struwe
\cite{SS95}. Uniqueness was proved only  under an additional
space-time integrability  of Strichartz type, which is
a crucial ingredient to the proof of the existence result. Indeed,
local existence can be proved using a fixed point argument in 
some Banach space which can be taken 
to be $B = C({\mathbb R};
H^1) \cap C^1({\mathbb R}; L^2) \cap L^5_{loc} ({\mathbb R};
L^{10})$. More recently, uniqueness was obtained under a different set
of conditions in \cite{Struwe}, using the energy inequality, but still with a space-time
integrability condition. One should also mention \cite{BG} where the
smooth solutions are proved to be globally in $L^5_t(L^{10}_x)$ and
stability under weak limits is proved.

In this paper, we intend to give a more physical condition 
which yields the uniqueness in the energy space. This condition 
can be easily understood in terms of finite speed of propagation.

We consider two solutions $u,v\in L^\infty(\dot H^1)\cap 
\dot{W}^{1,\infty}_t(L^2) $ to the 
wave equation   (\ref{eq:5d}),
with the same (real) initial  data $\phi_0\in \dot H^1$,
$\phi_1\in L^2$, namely
$$ u(t=0) = v(t=0) = \phi_0 \quad \quad  \partial_t u(t=0) = \partial_t v(t=0) =
\phi_1 ,$$
(note that the second condition, on $\partial_t u$ makes sense since a
solution which is in $ L^\infty(\dot H^1)\cap 
\dot{W}^{1,\infty}_t(L^2) $ is also in $C^1(\dot H^{-1})$).
 
 The two solutions $u$ and $v$  are supposed to be weak solutions, i.e. equation
(\ref{eq:5d}) holds in the  sense of distributions: for any
$\phi\in C^\infty_0([0,\infty);\R^3)$,
$$
\int_{\R\times \R^3} u \Box \phi+u^5\phi\, dxdt=\int_{\R^3} -\phi_1
(x)\phi(0,x)+\phi_0(x)\partial_t \phi(0,x)\,dx
$$
and the same equation holds for $v$.

\subsection{The local energy condition}
Let us state the local energy or the   finite speed of propagation condition. 
  Let $(t_0,x_0)$ be the vertex of   a backward cone $K$,
$K = \{ |x-x_0| = t_0 - t \} $ 
and $e(u)=|\partial u|^2/2+u^6/6$ be the energy density (here and
thereafter $\partial$ denotes the full space-time gradient). Then
we  assume that  for all $0\leq s \leq t \leq t_0$
\begin{eqnarray} 
  \label{eq:lei1}
  \int_{B(x_0,t_0-t)} e(u(t,x)) dx 
 \leq C    \int_{B(x_0,\alpha(t_0-s))} e(u(s,x)) dx.
\end{eqnarray}
where $C$ and $\alpha$ are some constants, $C \geq 1 $ and $\alpha \geq 1$.  

Similarly, consider the forward cone $K_1$ of vertex 
$(t_1,x_1)$, namely $K_1 = \{ |x-x_1| = t -t_1\} $, and let $t'\geq
t_1$, denote by $\partial_{K_1}$ the tangential derivatives, we assume that 
\begin{eqnarray} 
  \label{eq:lei}
 \frac 1 {\sqrt 2} \int_{t_1}^{t'} \int_{\partial 
 B(x_1,\tau-t_1)}
  \frac{|\partial_{K_1} u(\tau)|^2} 2+\frac{u(\tau)^6} 6  \quad   d\sigma d\tau 
 \\
 \leq  C  \int_{B(x_1,\alpha(t'-t_1))} e(u(t',x)) dx.   \nonumber 
\end{eqnarray}  
We insist on the fact that both constants $C$ and 
  $\alpha$   are supposed to be uniform with respect to the
vertex.

We point out that (\ref{eq:lei1})
  and  (\ref{eq:lei})  are  weak versions of the 
local energy  equality which is recalled in the next 
section. Indeed, for smooth solutions, one can prove that 
  (\ref{eq:lei1})
  and  (\ref{eq:lei}) hold  with $C=\alpha =1$. 
  We also notice that the left hand side of 
(\ref{eq:lei}) does not make sense (actually can be a priori infinite) 
if we only assume that $u \in  L^\infty(\dot H^1)\cap 
\dot{W}^{1,\infty}_t(L^2) $. 
Let us explain the meaning of  (\ref{eq:lei}). 
Let 
$\rho (x) \in C^\infty_0 (\R^3)$ be such that 
 $ \rho \geq 0$,  $\int \rho  = 1$ and define $\rho_n (x)  = n^3 
\rho(nx)$, then we define
  $u_n=
 u * \rho_n  $  a (space) regularization of $u$. Since 
$u \in \dot{W}^{1,\infty}_t(L^2)$, we deduce that $u_n$ is 
continuous in both space and time  variables. Condition  (\ref{eq:lei}) can 
be understood as 
\begin{eqnarray} 
  \label{eq:lei-reg}
 \limsup_{n\to \infty} \int_{t_1}^{t'} \int_{\partial 
 B(x_1,\tau-t_1)}
  \frac{|\partial_{K_1} u_n(\tau)|^2} 2+\frac{u_n(\tau)^6} 6 
  \quad   d\sigma d\tau   \\  
 \leq  C  \int_{B(x_1,\alpha(t'-t_1))} e(u(t',x)) dx.   \nonumber 
\end{eqnarray}

Let us prove that these conditions hold for any weak 
solution which also satisfies the local energy 
identity, namely 

\begin{equation} \label{energy-local}
\partial_t  e(u(t,x))  -  div  (\partial_t u \nabla u ) = 0.    
\end{equation} 

Let us prove that \eqref{eq:lei1} holds. We denote 
$M_s^t = \{ (\tau, x)  |  s<\tau<t , \,  |x-x_0| <  t_0 - \tau     \} $. 

Integrating \eqref{energy-local} over $M_s^t$, we formally get 
\eqref{eq:2lei1}.  
Let us prove this rigorously. Let 
$\rho (x) \in C^\infty_0 (\R^3)$ be such that 
 $ \rho \geq 0$,  $\int \rho  = 1$ and define $\rho_n (x)  = n^3 
\rho(nx)$. 
Hence    
\begin{equation} \label{energy-local-n} 
\partial_t  e(u(t,.)) * \rho_n    - 
 div  ( (\partial_t u \nabla u) *\rho_n ) = 0.    
\end{equation}
Using the fact that $ \partial_t u \nabla u  \in L^\infty_t (L^1_x) $, 
we deduce that $ (\partial_t u \nabla u) *\rho_n   \in L^\infty_t (C^\infty_x) $.
Hence, $  e(u(t,.)) * \rho_n   \in W^{1,\infty}_t (C^\infty_x) $. 
Integrating \eqref{energy-local-n} in $M_s^t$,  we get 

\begin{multline} 
  \label{eq:2lei1-reg}
  \int_{B(x_0,t_0-t)} e(u)* \rho_n \,dx+ 
 \frac 1 {\sqrt 2} \int_s^t \int_{\partial B(x_0,t_0-\tau)}
  [  e(u)-  (\partial_t u \nabla u) ] *\rho_n      \,   d\sigma
 d\tau \nonumber \\
      =    \int_{B(x_0,t_0-s)} e(u)* \rho_n   \,dx,  
\end{multline} 
 
Taking the limit when $n$ goes to infinity, we see that the 
fisrt and third terms converge to the corresponding terms in 
\eqref{eq:2lei1}.  For the second term, we rewrite 
$  e(u)-  (\partial_t u \nabla u) $ as 
$  e(u)-  (\partial_t u \nabla u) = 
  \frac{|\partial_K u |^2} 2+\frac{u ^6} 6   $. Then, 
using Jensen inequality, we deduce that 

\begin{equation*} 
  \frac{|\partial_{K} u_n(\tau)|^2} 2+\frac{u_n(\tau)^6} 6 
 \leq [ \frac{|\partial_{K} u (\tau)|^2} 2+\frac{u (\tau)^6} 6 ]* \rho_n.
\end{equation*}

Hence, 
\begin{multline} 
  \label{eq:2lei1-reg-limit}
  \int_{B(x_0,t_0-t)} e(u) \,dx+ \limsup_{n\to \infty}  
   \int_s^t \int_{\partial B(x_0,t_0-\tau)}
    \frac{|\partial_{K} u_n (\tau)|^2} 2+\frac{u_n (\tau)^6} 6   \,   d\sigma
 d\tau \nonumber \\
   \leq    \int_{B(x_0,t_0-s)} e(u)  \,dx.
\end{multline}

Arguing in the same way for the forward cone $K_1$, 
we deduce that  
\eqref{eq:lei-reg} holds with $C=\alpha = 1$.

\subsection{The main result}
We now state our main result.
\begin{theoreme} \label{th1} 
  Let $u$ be a weak solution to (\ref{eq:5d}) which satisfies
  (\ref{eq:lei1}) and  (\ref{eq:lei}).
 Then this solution is unique among all weak
  solutions satisfying (\ref{eq:lei1}) and  (\ref{eq:lei}).
\end{theoreme}
This unique solution is actually equal to the solution 
constructed in \cite{SS95}, but we will not use this fact in the proof, 
unlike for higher dimensions where a strong-weak uniqueness argument is 
used (\cite{fabwave} and remark at the end of the present paper), see also \cite{MN03} for a similar uniqueness 
result. 

  It does seem fairly reasonable for
 weak solutions to assume  that (\ref{eq:lei1})
  and  (\ref{eq:lei})  hold:
 certainly one is willing to have at least the weak energy
  inequality, namely $\int e(\phi)dx\leq \int e(\phi_0)dx$, and in
  light of the finite speed of propagation,both (\ref{eq:lei1}) and
 (\ref{eq:lei}) are not
  really stronger requirements. At any rate, control of the flux is already
  an essential tool in order to prove regularity for smooth data 
 (\cite{Gri,SSan}).
 
  A weak solution to (\ref{eq:5d}) satisfying in addition  (\ref{eq:lei})
  and  (\ref{eq:lei1})  can be considered as a {\em suitable 
  weak solution}. This is similar in spirit to the notion of suitable 
 weak solutions for the Navier-Stokes system introduced in 
  \cite{CKN83}. Indeed, both conditions are  local versions 
of the energy inequalities. 

\vspace{1mm} 
 
To prove  theorem \ref{th1},  we introduce a dual 
problem as was done in \cite{LM01cpde}. Then, we prove 
the existence of a smooth   solution to a regularized version of this 
 dual problem. This solution is used as a test 
function in the weak formulation. Passing to the limit, 
we deduce that $u=v$. 
 
In the next section, we recall the energy identities on  
 backward and forward  cones. In section 3, we give 
the proof  of theorem  \ref{th1}. We will start by a formal argument 
and then explain the regularization procedure.

\section{Finite speed of propagation} 

  Let us recall that a 
smooth solution of the wave equation (\ref{eq:5d}) satisfies
the following energy identity on each 
 backward
cone :  let again $(t_0,x_0)$ be the vertex of such a backward cone $K$, 
$K = \{ |x-x_0| = t_0 - t \} $ 
and $e(u)=|\partial u|^2/2+u^6/6$ be the energy density. Then
we have for all $s \leq t \leq t_0$
\begin{eqnarray} 
  \label{eq:2lei1}
  \,\,\,\int_{B(x_0,t_0-t)} e(u(t,x)) \,dx+ 
 \frac 1 {\sqrt 2} \int_s^t \int_{\partial B(x_0,t_0-\tau)}
  \frac{|\partial_K u(\tau)|^2} 2+\frac{u(\tau)^6} 6 \quad   d\sigma
 d\tau \nonumber \\
\quad \quad \quad \quad 
   =    \int_{B(x_0,t_0-s)} e(u(s,x)) \,dx,  
\end{eqnarray} 
where we recall that $\partial_K$ denotes the derivatives tangent to the backward
cone $K$.
 The second term on the left-hand side is usually referred to as
the (outgoing) flux through the cone $K$.

Moreover, the solution verifies the same inequality for
forward cones as well: specifically, consider the forward cone $K_1$ of vertex 
$(t_1,x_1)$, namely $K_1 = \{ |x-x_1| = t -t_1\} $, and let $t'\geq
t_1$, we have
\begin{eqnarray} 
  \label{eq:2lei}
 \frac 1 {\sqrt 2} \int_{t_1}^{t'} \int_{\partial 
 B(x_1,\tau-t_1)}
  \frac{|\partial_{K_1} u(\tau)|^2} 2+\frac{u(\tau)^6} 6 
  \quad   d\sigma d\tau 
 \\
  =    \int_{B(x_1,t'-t_1)} e(u(t',x)) \,dx.   \nonumber 
\end{eqnarray} 
The  left-hand side is usually referred to as
the (incoming) flux through the cone $K_1$.
\begin{rem} 
  Of course (\ref{eq:2lei}) is only a special case of the backward
  version of (\ref{eq:2lei1}): we could have stated an inequality
  between the two space-like surfaces $t=t'$ and $t=t''$ with $t_1\leq
  t''\leq t'$. Here we chose to take $t''=t_1$ as this is what will
  actually be needed later in the proof. 
\end{rem}

The conditions  (\ref{eq:lei})
  and  (\ref{eq:lei1})  which imply the uniqueness are 
 weaker versions of  (\ref{eq:2lei})
  and  (\ref{eq:2lei1}). Indeed, the equality is replaced 
by an inequality and we can even allow the presence of 
 fixed  constants $C$ and $\alpha$.

  Alternatively, one could rephrase both equalities in terms of only
  one equality, if one is willing to replace space balls by annuli
  (or even, say, domains with reasonably smooth boundaries). Then, if
  $\Sigma$ is the boundary of the backward domain of influence, one would
  ask the sum of the energy in our space domain at time $T$ and the
  outgoing flux through $\Sigma$ between times $T$ and $S\leq T$ to be
  equal to the energy at time $S$ in the space domain $\Sigma\cap
  \{t=S\}$. Such equalities and their weaker
  counterparts are a reasonable way to
  quantify the finite speed of propagation which one expects from any
  physically meaningful solutions to the equation.

\section{Proof of Theorem \ref{th1}}
\label{sec:sk}

Assume that $u$ and $v$ are two solutions of  (\ref{eq:5d}). Taking
$\phi$ to be an admissible test function $\phi\in
C^2_0([0,\infty),\R^3)$, we have
\begin{equation}
  \label{eq:weakuniq}
\int (u-v)\Box \phi+(u^5-v^5)\phi =0,  
\end{equation}
which can be rewritten as
\begin{equation}
  \label{eq:sk1}
  \int (u-v) (\Box \phi+(u^4+4u^3v+6u^2v^2+4uv^3+v^4)\phi)=0.
\end{equation}
We intend to solve the following (dual) problem: let $F\in
C_0^\infty((0,T)\times \R^3)$ and $\phi $ be the solution of the
following backward wave equation

\begin{equation}
  \label{eq:dual} \left\{ \begin{array}{l}
  \Box \phi+ V\phi=F, \\
  \phi(T)=\partial_t\phi(T)=0,  \end{array} \right.
\end{equation}
where we define $V=u^4+4u^3v+6u^2v^2+4uv^3+v^4$ and $T > 0$ is 
small enough, to be fixed later. Provided we solve
(\ref{eq:dual}) and prove that 
$\phi$ is regular enough to be used as a test function in (\ref{eq:weakuniq}),
 we will have uniqueness for our problem.  All is
required is for $\phi$ to be an admissible test function, in order
to justify the integration by parts. Actually, this will turn out 
to be untrue, but one may
still proceed using a  smoothing and a limiting procedure 
which will be explained later. 
 
\begin{proposition}\label{prop-inf}
  Provided $T$ is small enough, there exists a (compactly supported)
 smooth 
  solution $\phi_n$ to (a regularized version of)
 the dual problem (\ref{eq:dual}), such that
  $\phi_n$ is uniformly bounded in  $ L^\infty_{t,x}$.
\end{proposition}

\subsection{Formal proof}
Let us start by a formal proof. We will need a regularization of 
(\ref{eq:dual}) to make it rigorous. 
We denote   $K(z_0)$   the forward cone with vertex $z_0=(t_0,x_0)$ 
and time $t \leq T$
i.e. $K(z_0)=\{(t,x) | \  |x-x_0|=t-t_0, \ t_0 \leq t \leq T \}$. 
 Then, the solution of (\ref{eq:dual}) is 
given by, taking advantage of the explicit space
representation of the fundamental solution to the 3D wave equation,
\begin{equation}
  \label{eq:integrale}
  \phi(t_0,x_0)=\int_{K(z_0)} 
 \frac{F(z)-V\phi(z)}{|z-z_0|}d\sigma(z),
\end{equation}
with $z = (t,x)$ and 
where $\sigma$ is the surface measure on forward cones. Then, we proceed as
J\"orgens (\cite{Jo}), with
\begin{equation} \label{phiLinfinie} 
\|\phi\|_{L^\infty((0,T)  \times \R^3)} 
 \leq C(F)+\|\phi\|_{L^\infty((0,T)  \times \R^3)}  \sup_{z_0} \int_{K(z_0)}
\frac{|V(z)|}{|z-z_0|}d\sigma(z),
 \end{equation}
and as $|V|\lesssim u^4+v^4,
$
we use
$$
\int_K \frac{|u|^4}{|z-z_0|}d\sigma(z) \sim
\int_{B(0,T-t_0)}\frac{|w(y)|^4}{|y|} dy,
$$ 
where $w(y)=u(t_0+|y|,x_0+y)$. This in turn yields
\begin{equation}
  \label{eq:bla1}
\int_K \frac{|u|^4}{|z-z_0|}d\sigma(z) \lesssim
\int_{B(0,T-t_0)}\frac{|w(y)|^2}{|y|^2} dy+ \int_{B(0,T-t_0)}|w|^6.  
\end{equation}
By an appropriate local version of Hardy's inequality (see
e.g. \cite{SSbook}), the first term in (\ref{eq:bla1}) is controlled:
\begin{equation}
  \label{eq:bla2}
  \int_{B(0,T-t_0)}\frac{|w(y)|^2}{|y|^2} dy  \lesssim 
\int_{B(0,T-t_0)}|\nabla_y w|^2 dy+ (\int_{B(0,T-t_0)}|w(y)|^6
dy)^{2/6}.
\end{equation}
We then recognize the flux,
$$
\frac{1}{2}\int_{B(0,t_0)}|\nabla_y w|^2+
\frac{1}{6}\int_{B(0,t_0)}|w|^6\, dy = \mathrm{flux},
$$
as
$$
|\nabla_y w|^2=|\nabla u-\frac y {|y|} \partial_t u|^2,
$$
and recall
$$
\mathrm{flux}=\int_{K}\frac{1}{2}\left|\nabla u-\frac y {|y|}
  \partial_t u\right |^2+\frac{|u|^6}{6} d\sigma.
$$
Hence (\ref{eq:bla2}) becomes
\begin{equation}
  \label{eq:bla3}
  \int_K \frac{|u|^4}{|z-z_0|}d\sigma(z) \lesssim
  \mathrm{flux}+\mathrm{flux}^{\frac 1 3}.
\end{equation}
By choosing $T$ small enough, we can make the local energy
$\int_{B(x_0,\alpha T)} e(u(T,x))dx$ smaller than a fixed constant
$\varepsilon_0$, uniformly in $x_0$: we simply use the energy
inequality (\ref{eq:lei1}), fixing $T$ such that
$\int_{B(x_0,\alpha(\alpha+1)  T)}e(u(0,x))dx$ is (uniformly) small enough, which in
turn is a trivial consequence of the initial  data  being in  $\dot H^1\times
L^2$. Then we deduce that the flux through the forward cone which is
needed in the construction of $\phi$ can be made  smaller than $1/2$ by using
(\ref{eq:lei}) and choosing $\varepsilon_0$ such that $C \varepsilon_0 =\frac12  $.
Next,  we  
can perform a contraction argument in $L^\infty_{t,x}$ to obtain
$\phi$. 
\begin{rem}
  Note that the whole argument is local in space-time. Hence, the assumptions on the data could be relaxed to $\dot
  H^1_{\text{loc}}\times L^2_{\text{loc}}$, and one could consider local in time weak solutions. We elected to
  keep $\dot H^1\times L^2$ data and global in time solutions for simplicity.
\end{rem}
\subsection{Rigorous proof} 
Let us  explain the regularization procedure which yields 
a rigorous proof of (\ref{phiLinfinie} ) and  the proposition. 
 Recall that $\rho (x) \in C^\infty_0 (\R^3)$ is such that 
 $ \rho \geq 0$,  $\int \rho  = 1$ and $\rho_n (x)  = n^3 
\rho(nx)$, then we define $u_n=u*\rho_n$, $v_n=v*\rho_n$ and $V_n = 
u_n^4+4u_n^3v_n+6u_n^2v_n^2+4u_nv_n^3+v_n^4$. We intend to solve
\begin{equation} 
  \label{eq:dual2} \left\{ \begin{array}{l}
  \Box \phi_n+ V_n  \phi_n =F, \\
  \phi_n(T)=\partial_t\phi_n(T)=0,  \end{array} \right.
\end{equation}
by a fixed point argument. Considering
\begin{equation} 
  \label{eq:dual2bis} \left\{ \begin{array}{l}
  \Box \psi+ V_n  \tilde\psi =F, \\
  \psi(T)=\partial_t\psi(T)=0,  \end{array} \right.
\end{equation}
for smooth $\psi$ and $\tilde\psi$, we have
\begin{equation}
  \label{eq:integrale-reg}
  \psi(t_0,x_0)=\int_{K(z_0)} 
 \frac{F(z)- (V_n\tilde\psi) (z)}{|z-z_0|}d\sigma(z),
\end{equation}
from which we infer that 
\begin{equation} \label{phiLinfinie-reg} 
\|\psi\|_{L^\infty((0,T)  \times \R^3)} 
 \leq C(F)+\|\tilde \psi\|_{L^\infty((0,T)  \times \R^3)}  \sup_{z_0} \int_{K(z_0)}
\frac{|V_n|}{|z-z_0|}d\sigma(z).
 \end{equation}
Now, we can proceed as in the formal proof and choose $T$ small enough so 
that 
$$  \sup_{n}  \sup_{z_0} \int_{K(z_0)}
\frac{|V_n|}{|z-z_0|}d\sigma(z)   \leq \frac12    . $$

Notice that given we are solving a linear problem, estimating
$\psi$ or $\psi- \phi$ is identical,  where $\phi$ solves 
\begin{equation} 
  \label{eq:dual2bis2} \left\{ \begin{array}{l}
  \Box \phi+ V_n  \tilde\phi =F, \\
  \phi(T)=\partial_t\phi(T)=0.   \end{array} \right. 
\end{equation} 
Hence we deduce from the previous computations  that 
$$ 
\|\psi - \phi\|_{L^\infty_{t,x}}\leq \frac 1 2
\|\tilde \psi-\tilde\phi\|_{L^\infty_{t,x}}.
$$ 
This estimate 
allows  a fixed point argument in  $C^0$ to be carried out. Therefore we
 have constructed  a solution $\phi_n$ to the equation (\ref{eq:dual2}). Moreover,
 we recover an estimate on $\|\phi_n\|_{L^\infty_{t,x}}$ which is
 uniform with respect to $n$, thanks to
 (\ref{eq:lei}). Furthermore, $\phi_n$ is smooth, as the regularity
 can as usual be carried along the iterates which yield $\phi_n$: for
 any derivative $\partial$, we have
$$ 
\|\partial \psi - \partial \phi\|_{L^\infty_{t,x}}\leq \frac 1 2
\|\partial \tilde \psi-\partial \tilde\phi\|_{L^\infty_{t,x}}+\|\tilde
\psi- \tilde\phi\|_{L^\infty_{t,x}} C(\partial V_n).
$$ 
 We
 do not get good control of norms, as they involve derivative of
 $V_n$, but we will not need it. Moreover, $\phi_n$ is compactly supported, by finite
 speed of propagation (again, all iterates are in a uniform way).
 This ends the proof of   proposition \ref{prop-inf}.

We now return to the proof of Theorem \ref{th1} and explain why the
smoothing procedure which yields $\phi_n$ still allows for the
heuristic argument to be carried out. In fact, we use $\phi_n$ as a
test function: for all $n$, we have 
\begin{equation}
  \label{eq:sk2}
  \int (u-v) ( \Box \phi_n  +V  \phi_n  )=0.
\end{equation} 
This translates into
\begin{equation}
  \label{eq:sk22} 
  \int (u-v) ( F  +(V-V_n) \phi_n  )=0.
\end{equation} 
We know that $u-v\in L^\infty_t(L^6)$, and that  $  V_n  $
converges to $V$ strongly in
$L^\infty_{t,\mathrm{loc}}L^{3/2}$. Hence,  $(V-V_n)\phi_n$   converges  toward zero
in $L^1_{t,\mathrm{loc}}L^{6/5}_{x,\mathrm{loc}}$, given that
$\phi_n$ is uniformly in $L^\infty_{t,x} $ and is compactly supported;
this ultimately gives the desired equality:
$$ 
 \int (u-v) F=0, 
$$ 
from which we deduce that $u=v$ on the interval $[0,T]$. Now,   we can
argue by contradiction, choosing  the initial time to be   $t_0$
where 
$$  t_0 = \hbox{inf}  \{ t, \,  t \geq 0, \,  u(t) \neq  v(t)      \}     .  $$
Using the fact that $u$ and $v$ are continuous in time with values in 
$ \dot H^{-1}   $, we deduce that   $u(t_0)=v(t_0)$.  
 Then we have $u=v$ on some interval $[t_0,
t_0+\eta]$, which proves that no such $t_0$ exists. Hence, we deduce  that $u=v$ on $[0, \infty)$
which achieves the proof of the main theorem.

Finally, we make some comments on the case  $n\geq 4$.
 In higher dimensions, one cannot rely on J\"orgens
  estimate. However, uniqueness was proven under the assumption
  $\phi\in C_t(\dot H^1)\cap C^1_t(L^2)$ in \cite{fabwave}, for both
  focusing and defocusing critical wave equation, with $n\geq 4$. In
  the defocusing case, assuming only the local energy identity
  (\ref{eq:lei1}), one can easily get rid of the continuity in time and
  obtain uniqueness as in Theorem 1. We refer the interested reader to
  \cite{fn2} for further discussions in a similar (albeit more
  complicated) setting.

\section{Acknowledgments}
The authors would like to thank Jalal Shatah for many discussions 
about this work. The first author was partially supported by an 
NSF grant and by an Alfred Sloan Fellowship. Part of this work 
was done while the second author was visiting the Courant 
Institute, which he would like to thank for its hospitality. 


 \end{document}